\input amstex
\documentstyle{amsppt}

\magnification=\magstep 1

\NoBlackBoxes
%\hsize=6 true in 
%\hoffset=.25 true in

%%%%%%%%%%%%MACROS%%%%%%%%%%%%

\define\A{\Cal A}
\define\Ah{A_\hbar}
\redefine\B{\Cal B}
\define\Bh{B_\hbar}

\define\C{\Bbb C}
\define\CC{\Cal C}
\define\E{\Cal E}
\define\F{\Cal F}
\define\h{\hbar}
\redefine\H{\Cal H}
\define\I{\Cal I}
\define\K{\Cal K}

\define\N{\Bbb N}
\define\R{\Bbb R}

\define\Z{\Bbb Z}

\define\dbar{\bar{\partial}}

\define\id{\operatorname{id}}
\define\ideal{\triangleleft}
\define\Interior{\operatorname{Interior}}
\define\LUB{\operatorname{LUB}}
\define\>{\rangle}
\define\<{\langle}
\define\cospec{\operatorname{cospec}}
\define\spec{\operatorname{spec}}
\define\supp{\operatorname{supp}}
\define\rank{\operatorname{rank}}
\define\tensor{\otimes}
\define\tr{\operatorname{trace}}

\define\lcl{[\![}
\define\rcl{]\!]}
\define\({(\!(}
\define\){)\!)}

%%%%%%%%%%%END MACROS%%%%%%%%%%
\topmatter 
\title Asymptotic Spectral Measures, Quantum Mechanics, and $E$-theory \endtitle
\rightheadtext{Asymptotic Spectral Measures, Quantization, and E-theory}
\author Diane Martinez and Jody Trout \endauthor
\affil Dartmouth College \endaffil
\thanks Research of the second author was supported by NSF grant
DMS-0071120\endthanks
\address \newline 6188 Bradley Hall \newline
Dartmouth College \newline
Hanover, NH 03755 \newline
diane.martinez\@dartmouth.edu \newline
jody.trout\@dartmouth.edu \endaddress

\abstract 
We study the relationship between POV-measures in quantum theory and 
asymptotic morphisms in the operator algebra $E$-theory
of Connes-Higson. This is done by introducing the theory of
``asymptotic'' PV-measures
and their integral correspondence with positive asymptotic morphisms on locally
compact spaces. Examples and applications involving
various aspects of quantum physics, including quantum noise models, 
semiclassical limits, strong deformation quantizations,
and pure half-spin particles, are also discussed.
\endabstract
\endtopmatter

%%%%%%%%%%%%%%%%%%%%%%%%%%%%%%%%%%%%%%%%%%%%%%%%%%%%%%%%%%%%%%%%%%%%%%%%%%%%%%%%
%%%%%%%%%%%%%%%%%%%%%%%%%%%%%%%%%%%%%%%%%%
\head 1. Introduction \endhead
%%%%%%%%%%%%%%%%%%%%%%%%%%%%%%%%%%%%%%%%%%%%%%%%%%%%%%%%%%%%%%%%%%%%%%%%%%%%%%%%
%%%%%%%%%%%%%%%%%%%%%%%%%%%%%%%%%%%%%%%%%%

In the Hilbert space formulation of quantum mechanics by von Neumann \cite{VN},
an observable is modeled as a self-adjoint operator on the Hilbert 
space of states of the quantum system. The Spectral Theorem relates this
theoretical view of a quantum observable to the more 
operational one of a projection-valued measure (PVM or spectral measure) which
determines the probability distribution of the experimentally measurable values
of the
observable. To solve
foundational problems with the concept of measurement and to better analyze
unsharp results in experiments, this view  was
generalized to
include positive operator-valued measures (POVMs).  Since the work of Jauch and
Piron \cite{JP}, 
POV-measures have played an ever increasing role in both the foundations and
operational aspects of quantum physics  \cite{BGL,S}. 
See the Appendix for a quick review of POVMs, their use in quantum
mechanics, and relation to the Spectral Theorem.

In this paper, we study the relationship between POVMs and asymptotic morphisms
in the operator algebra $E$-theory of 
Connes and Higson \cite{CH}, which has already found many applications in
mathematics  \cite{Bl,GHT,H,Tr}, most notably to
classification problems in operator $K$-theory, index theory, representation
theory,
geometry, and topology.  The basic ingredients of $E$-theory
are asymptotic morphisms, which are given by continuous families of functions
$\{Q_\hbar\}_{\hbar >0} : \A \to \B$ from a
$C^*$-algebra $\A$ to a $C^*$-algebra $\B$ that satisfy the axioms of a
$*$-homomorphism in the limit as the parameter $\hbar$ tends to $0$. 
Asymptotic multiplicativity is a modern version of the Bohr-von Neumann
correspondence principle \cite{L} from quantization theory: For all $f, g \in
\A,$
$$
Q_\hbar(fg) - Q_\hbar(f)Q_\hbar(g) \to 0 \text{ as } \hbar \to 0.
$$
It is then no surprise that quantization schemes may naturally define
asymptotic  morphisms, say, from the $C^*$-algebra $\A$ of classical observables
to the
$C^*$-algebra $\B$ of quantum observables. 
Hence, such quantizations can give cycles in the abelian group $E(\A,\B)$, which was
defined by Connes and Higson as a certain matrix-stable
homotopy group of asymptotic morphisms from $\A$ to $\B$.
For example, Guentner \cite{G1}
showed that Wick quantization on the Fock space $\F$ of $\C$ defines a
positive
asymptotic morphism $\{Q^W_\hbar\} : C_0(\C) \to \K(\F)$, whose $E$-theory
class is equal to the class of the $\dbar$-operator $\lcl \dbar \rcl = \lcl
Q^W_\hbar \rcl \in E(C_0(\C), \C)$.
(We will discuss Guentner's work in our context in Example 5.5.)
See the papers \cite{N1,N2,Ro} and the books \cite{C,GVF} for more on the
connections
between operator algebra $K$-theory, $E$-theory, and quantization.  

We show that there is a fundamental quantum-$E$-theory relationship by
introducing the concept of an {\it asymptotic spectral measure} (ASM or
asymptotic PVM)
$$\{A_\hbar\}_{\hbar > 0} : \Sigma \to \B(\H)$$
associated to a  measurable space $(X, \Sigma)$. (See Definition 3.1.) Roughly,
this is a continuous
family of POV-measures which are ``asymptotically'' projective (or
quasiprojective)
as $\hbar$ tends to $0$:
$$ \Ah(\Delta)^2 - \Ah(\Delta) \to 0 \text{ as  } \hbar \to 0$$
for certain measurable sets $\Delta \in \Sigma$.

Let $X$ be a locally compact space with Borel $\sigma$-algebra $\Sigma_X$
and let $C_0(X)$ denote the $C^*$-algebra
of continuous functions vanishing at infinity on $X$. One of our main results is
an ``asymptotic'' Riesz representation theorem (Theorem 4.2) which gives a
bijective correspondence between certain 
positive asymptotic morphisms $$\{Q_\hbar\} : C_0(X) \to \B$$ and Borel
asymptotic spectral measures 
$$\{\Ah\} :  (\Sigma_X,\CC_X ) \to  (\B(\H), \B)$$
where $\CC_X$ denotes the open subsets of $X$ with compact closure and $\B$ is a
hereditary $C^*$-subalgebra of $\B(\H)$. This correspondence is
given by operator integration
$$Q_\hbar(f) = \int_X f(x) \, d\Ah(x).$$
The associated asymptotic morphism $\{Q_\hbar\} : C_0(X) \to \B$ then allows
one to define an $E$-theory invariant for the asymptotic spectral measure
$\{\Ah\}$,
$$ \lcl \Ah \rcl =_{\text{def}} \lcl Q_\hbar \rcl \in E(C_0(X), \B)  \cong
E_0(X; \B),$$
in the $E$-homology group of $X$ with coefficients in $\B$. 

It has been well-established that operator $K$-theory and the dual $K$-homology
groups provide suitable receptacles for
invariants of quantum systems, such as chiral anomalies in quantum field theory
\cite{N} and, more recently, as $D$-brane charges in string theory and $M$-theory \cite{P,W}.
Since $E$-theory subsumes both $K$-theory and $K$-homology \cite{Bl},
it is reasonable to assume that $E$-theory elements of quantizations
and asymptotic spectral measures may provide interesting topological invariants
of 
the associated quantum systems. Although in this paper we will be more concerned
with asymptotic morphisms
and their relation to POV-measures than computing $E$-theory elements (but see
Example 5.5), a long-range
goal of this research project is to develop an $E$-theoretic calculus for
computing these invariants directly from the asymptotic measure-theoretic data, e.g.,
by developing the appropriate notions of homotopy and suspension for ASMs, thus 
bypassing the technical functional-analytic aspects of asymptotic morphisms.

Another benefit of using this asymptotic measure-theoretic approach is operational
in nature. Experimental data from position and momentum measurements on
an elementary quantum system (via visibility data from interference experiments)
is collected which is then used to
construct the associated POVM. This method
\cite{S} is based on using frame manuals for the instrument state space and Sakai operators associated
to localization operators on rectangles in the classical phase space $X$. The POVM $\{\Ah\}$
depends on Planck's constant, of course, and generally satisfies the (unsharp) separation property
$$\Ah(\Delta_1 \cap \Delta_2) \neq \Ah(\Delta_1) \Ah(\Delta_2).$$
However, if letting $\hbar \to 0$ one then obtains an ASM, which is equivalent to
$$\lim_{\hbar \to 0} \left( \Ah(\Delta_1 \cap \Delta_2) - \Ah(\Delta_1) \Ah(\Delta_2) \right) = 0,$$
then one can directly associate an $E$-homological invariant $\lcl \Ah \rcl \in E_0(X; \B)$
to the quantum system under experimental study using our theory.

The outline of this paper is as follows. In Section 2 we discuss positive
asymptotic morphisms associated
to hereditary and nuclear $C^*$-algebras. The basic definitions and properties
of asymptotic spectral measures 
are developed in Section 3. Asymptotic Riesz representation theorems and some of
their consequences
are proven in Section 4. Examples and applications of ASMs associated to various
aspects of quantum physics are discussed in Section 5, e.g.,
constructing ASMs from PVMs by quantum noise models, quasiprojectors and
semiclassical limits, unsharp spin measurements of spin-$\frac12$ particles 
(including an example from quantum cryptography), strong deformation
quantizations, and Wick quantization on bosonic Fock space.

The authors would like to thank Navin Khaneja, Iain Raeburn, and Dana Williams
for helpful conversations. Also, we would like to thank the referee for helpful
comments. See Beggs \cite{B} for a related method of obtaining asymptotic
morphisms
by an integration technique involving spectral measures.

%%%%%%%%%%%%%%%%%%%%%%%%%%%%%%%%%%%%%%%%%%%%%%%%%%%%%%%%%%%%%%%%%%%%%%%%%%%%%%%%
%%%%%%%%%%%%%%%%%%%%%%%%%%%%%%%%%%%%%%%%%%
\head 2. Positive Asymptotic Morphisms and Hereditary $C^*$-subalgebras \endhead
%%%%%%%%%%%%%%%%%%%%%%%%%%%%%%%%%%%%%%%%%%%%%%%%%%%%%%%%%%%%%%%%%%%%%%%%%%%%%%%%
%%%%%%%%%%%%%%%%%%%%%%%%%%%%%%%%%%%%%%%%%%

Let $\A$ and $\B$ be $C^*$-algebras. Recall that a linear map $Q : \A
\to \B$
is called positive \cite{M} if $Q(f) \geq 0$ for all $f \geq 0$ in $\A$. It is
called
completely 
positive if every inflation  to $n \times n$ matrices $M_n(Q) : M_n(\A) \to
M_n(\B)$ is also 
positive. Every $*$-homomorphism from $\A$ to $\B$ is clearly completely
positive. The following definition
interpolates between (completely) positive linear maps and $*$-homomorphisms.

\definition{Definition 2.1} A ({\it completely}) {\it positive asymptotic
morphism} from $\A$ to $\B$ is a
family of maps
$$\{Q_\hbar\}_{\hbar \in (0,1]} : \A \to \B$$
parameterized by  $\hbar \in (0,1]$ such that the following conditions hold:
\roster
\item"a.)" Each $Q_\hbar$ is a (completely) positive linear map; 
\item"b.)" The map $(0,1] \to \B: \hbar \to Q_\hbar(f)$ is continuous for each
$f \in \A$;
\item"c.)" For all $f, g \in \A$ we have $$\lim_{\hbar \to 0} \| Q_\hbar(fg) -
Q_\hbar(f)Q_\hbar(g)\| = 0.$$
\endroster
\enddefinition For the basic theory of asymptotic morphisms see the books
\cite{GHT,C,Bl} and papers \cite{CH,G2}. For the importance
of positive asymptotic morphisms to $C^*$-algebra $K$-theory see \cite{HLT}.
Note that any $*$-homomorphism
$Q : \A \to \B$ determines the constant completely positive asymptotic morphism
$\{Q_\hbar\} : \A \to \B$
defined by $Q_\hbar = Q$ for all $\hbar > 0$. Also, it follows that for any
$f\in \A$,
a mild boundedness condition \cite{Bl} always holds,
$$\limsup_{\hbar \to 0} \| Q_\hbar(f) \| \leq \|f\|.$$

\remark{Remark} In the $E$-theory literature, asymptotic morphisms are usually
parameterized by $t \in [1, \infty)$. We chose to
use the equivalent parameterization $\hbar = 1/t \in (0,1]$ to make the
connections to quantum physics more transparent. Note that
other authors have used different parameter spaces, including discrete ones
\cite{L2,Th}. The results in this paper translate verbatim
to these parameter spaces, and condition (b.) is obviously irrelevant in the
discrete
case.
\endremark

\definition{Definition 2.2} Two asymptotic morphisms $\{Q_\hbar\}, \{Q_\hbar'\}
: \A \to \B$ are called {\it equivalent}
if for all $f \in \A$ we have that $\lim_{\hbar \to 0} \| Q_\hbar(f) -
Q_\hbar'(f) \| = 0.$ We will let $\lcl \A, \B\rcl_{a(cp)}$ denote the collection
of all asymptotic equivalence
classes of (completely positive) asymptotic morphisms from $\A$ to $\B$.
\enddefinition

A $C^*$-algebra $\A$ is called nuclear \cite{M} if the identity map $\id : \A
\to \A$ can be approximated pointwise in norm
by completely positive finite rank contractions. This is equivalent to the
condition that there is a 
unique $C^*$-tensor product $\A \tensor \B$ for any $C^*$-algebra $\B$. If $\H$
is a separable Hilbert space, the
$C^*$-algebra $\K(\H)$ of compact operators on $\H$ is nuclear. If $X$ is a
locally compact space, then the $C^*$-algebra
$C_0(X)$ of continuous complex-valued functions on $X$ vanishing at infinity is
also nuclear.

If $\A \cong C(X)$ is unital and commutative, then every positive
linear map $Q : \A \to \B$ is completely positive by Stinespring's Theorem.
The following result is a consequence of the completely positive lifting theorem
of Choi and Effros \cite{CE} for nuclear $C^*$-algebras. 
(See also 25.1.5 of Blackadar \cite{Bl} for a discussion.)

\proclaim{Lemma 2.3} Let $\A$ be a nuclear $C^*$-algebra. Every asymptotic
morphism from $\A$ to any $C^*$-algebra $\B$ is
equivalent to a completely positive asymptotic morphism. That is, there is a
bijection of sets
$\lcl \A, \B \rcl_a \cong \lcl \A, \B \rcl_{acp}.$
\endproclaim

\definition{Definition 2.4} Let $\A_1 \subset \A$ and $\B_1 \subset \B$ be 
subalgebras of the $C^*$-algebras $\A$ and $\B$. If $Q : \A \to \B$ is a
linear map such that $Q(\A_1) \subset \B_1$, we will denote this
by
$$Q :   (\A,\A_1) \to  (\B,\B_1).$$
The notation $\{Q_\hbar\} : (\A,\A_1) \to  (\B,\B_1)$ then has the
obvious meaning.
\enddefinition

\proclaim{Lemma 2.5} Let $\A_1 \subset \A$ and $\B_1 \subset \B$ be non-closed
$*$-subalgebras.
Every positive linear map $Q : (\A,\A_1) \to  (\B,\B_1)$
also satisfies 
$$Q : (\A,\overline {\A_1}) \to  (\B, \overline{\B_1}),$$
where $\overline{\A_1}$ denotes the closure of $\A_1 \subset \A$ (similarly for
$\overline{\B_1}$).
\endproclaim 

\demo{Proof} Follows from the fact that a positive linear map is automatically
norm bounded.
\qed
\enddemo

Let $\A$ be a $*$-subalgebra of a $C^*$-algebra $\B$. Recall that $\A$ is said
to be hereditary \cite{M} if $0 \leq b \leq a$ and
$a \in \A$ implies that $b \in \A$. Every (closed two-sided $*$-invariant) ideal
in a $C^*$-algebra is a hereditary $*$-subalgebra. 
In particular, if $\H$ is a Hilbert space, the ideal of compact operators
$\K(\H)$ is a hereditary $C^*$-subalgebra of the
$C^*$-algebra of bounded operators $\B(\H)$. An important (non-closed)
hereditary $*$-subalgebra
for quantum theory is the (non-closed) ideal $\B_1(\H) \subset \K(\H)$ of
trace-class
operators: $$\B_1(\H) = \{ \rho \in \K(\H) : \tr |\rho| < \infty \}.$$
We then have that $\K(\H)^* = \B_1(\H)$ by the dual pairing $\rho(T) = \tr(\rho
T)$, where $\rho \in \B_1(\H)$ and $T \in \K(\H).$

If $X$ is a locally compact space, then the ideal $C_0(X)$ is a hereditary
$C^*$-subalgebra of the 
$C^*$-algebra $C_b(X)$ of  continuous bounded complex-valued functions on $X$.
Also, the (non-closed) ideal $C_c(X)$ 
of compactly supported functions is a (non-closed) hereditary $*$-subalgebra of
$C_b(X)$. However, in general, $C_\delta(X)$, for $\delta = $ c, $0$, b,
is not a hereditary subalgebra of the $C^*$-algebra $B_b(X)$ of bounded Borel
functions on $X$.

%%%%%%%%%%%%%%%%%%%%%%%%%%%%%%%%%%%%%%%%%%%%%%%%%%%%%%%%%%%%%%%%%%%%%%%%%%%%%%%%
%%%%%%%%%%%%%%%%%%%%%%%%%%%%%%%%%%%%%%%%%%
\head 3. Asymptotic Spectral Measures \endhead
%%%%%%%%%%%%%%%%%%%%%%%%%%%%%%%%%%%%%%%%%%%%%%%%%%%%%%%%%%%%%%%%%%%%%%%%%%%%%%%%
%%%%%%%%%%%%%%%%%%%%%%%%%%%%%%%%%%%%%%%%%%

In this section we give the basic definitions and properties of asymptotic
spectral measures.  See the Appendix for a review of POV and spectral measures.
Let $(X, \Sigma)$ be a measurable space and $\H$ a separable Hilbert space. 
Let $\E \subset \Sigma$ denote a fixed collection of measurable subsets.

\definition{Definition 3.1}
A {\it asymptotic spectral  measure} (ASM) on $(X,\Sigma, \E)$ is a family of
maps 
$$\{\Ah\}_{\hbar \in (0,1]} : \Sigma \to \B(\H)$$
parameterized by $\hbar \in (0,1]$ such that the following conditions hold:
\roster
\item"a.)" Each $\Ah$ is a POVM on $(X,\Sigma)$ with $\limsup_{\hbar \to 0}
\|\Ah(X)\| \leq 1$;
\item"b.)" The map $(0,1] \to \B(\H) : \hbar \to \Ah(\Delta)$ is continuous for
each $\Delta \in \E$;
\item"c.)" For each $\Delta_1$, $\Delta_2 \in \E$ we have that
$$\lim_{\hbar \to 0} \|  \Ah(\Delta_1 \cap \Delta_2) - \Ah(\Delta_1)
\Ah(\Delta_2) \| = 0.$$
\endroster
The triple $(X, \Sigma, \E)$ will be called an {\it asymptotic measure space}.
The family $\E$
will be called the {\it asymptotic carrier} of $\{\Ah\}$.
Condition (c.) will be called {\it asymptotic projectivity} (or {\it
quasiprojectivity}) and generalizes the projectivity condition (A.1) of 
a spectral measure. It is motivated by the quantum theory notion of
quasiprojectors, as discussed in Example 5.2. If $\E = \Sigma$ then we will call
$\{\Ah\}$ a {\it full}
ASM on $(X, \Sigma)$. If each $\Ah$ is normalized, i.e., $\Ah(X) =I_\H$,
then we will say that $\{\Ah\}$ is {\it normalized}. The mild
boundedness condition in (a.) is then redundant. (Also see the remark after
Definition 2.1.)
\enddefinition 

A spectral (PV) measure $E : \Sigma \to \B(\H)$ determines a ``constant'' full
asymptotic spectral measure $\{\Ah\}$ by the assignment $\Ah = E$
for all $\hbar$. Also, any continuous family $\{E_\hbar\}$ of spectral measures
(in the sense of (b.)) determine an ASM on $(X,\Sigma, \E)$.
See \cite{CHM} for an application of smooth families of spectral measures to the
Quantum Hall Effect. 

\definition{Definition 3.2} Two asymptotic spectral measures $\{\Ah\}, \{\Bh\} :
\Sigma \to \B(\H)$
on $(X, \Sigma, \E)$ are said to be
(asymptotically) {\it equivalent} if for each measurable set $\Delta \in \E$,
$$\lim_{\hbar \to 0} \| \Ah(\Delta) - \Bh(\Delta) \| = 0.$$
This will be denoted $\{\Ah\} \sim_\E \{\Bh\}$. If this holds for $\E = \Sigma$
we will
call them {\it fully equivalent}.
\enddefinition

From now on, we let $X$ denote a locally compact Hausdorff topological space
with
Borel $\sigma$-algebra $\Sigma_X$. We will assume that $\E = \CC_X$ denotes the
collection of all open
subsets of $X$ with compact closure, i.e., the pre-compact open subsets.

\definition{Definition 3.3} Let $\B \subset \B(\H)$ be a hereditary
$*$-subalgebra.
A Borel POV-measure $A : \Sigma_X \to \B(\H)$
will be called {\it locally $\B$-valued} if $A(U) \in \B$ for all pre-compact
open
subsets
$U \in \CC_X$ and this will be denoted by
$$A : (\Sigma_X, \CC_X) \to (\B(\H), \B).$$
A family of Borel POV-measures $\{\Ah\}$ on $X$ will be called
{\it locally $\B$-valued} if each
POVM $\Ah$ is locally $\B$-valued and will be denoted
$\{\Ah\} :   (\Sigma_X,\CC_X) \to  (\B(\H),\B).$
We will use the term {\it locally compact-valued} for locally $\K(\H)$-valued.
If $\B =
\B_1(\H) \subset \K(\H)$ is the trace-class operators,
then we will say that $\{\Ah\}$ has {\it locally compact trace}.
\enddefinition

We will let $\( X, \B\)$ denote the set of all equivalence classes 
of locally $\B$-valued Borel asymptotic spectral measures on $(X,\Sigma_X,
\CC_X)$. The equivalence class
of $\{\Ah\}$ will be denoted $\(\Ah\) \in \(X, \B\)$.

Given a Borel POV-measure $A$ on $X$, the {\it cospectrum} of $A$ is defined as
the set
$$\cospec(A) = \bigcup \{U \subset X: U \text{ is open and } A(U) = 0\}.$$
The {\it spectrum} of $A$ is the complement $\spec(A) = X \backslash
\cospec(A)$. 
The following definition is adapted from Berberian \cite{Be}

\definition{Definition 3.4} A POVM $A$ on $X$ will be said to have {\it compact
support} if the spectrum of $A$ is a compact
subset of $X$. An ASM $\{\Ah\}$ on $X$ will be said
to have {\it compact support} if there is 
a compact subset $K$ of $X$ such that $\spec(\Ah) \subset K$ for all $\hbar >
0$.
\enddefinition

The relationship among these compactness notions is contained in the
following.

\proclaim{Proposition 3.5} Let $X$ be second countable. Let $A$
be a Borel POVM on $X$ with compact support. Let $\B$ be the 
hereditary subalgebra of $\B(\H)$ generated by $A(\spec(A))$. Then $A$ is a
locally $\B$-valued POVM, i.e., $A : (\Sigma_X,\CC_X) \to  (\B(\H),\B).$
\endproclaim

\demo{Proof}  Since $X$ is second countable, the $\sigma$-algebra $\B_X$ 
of Baire subsets equals the Borel $\sigma$-algebra $\B_X = \Sigma_X$. Thus, by
Theorem 23 \cite{Be}
$A(\cospec(A)) = 0$. Let $U \in \CC_X$ be a pre-compact open subset of $X$.
We then have that $$0 \leq A(U \cap \cospec(A)) \leq A(\cospec(A)) = 0$$ and
since $X$ is the disjoint union $X = \spec(A) \sqcup \cospec(A)$,
$$0 \leq A(U) =  A(U \cap \spec(A)) \leq A(\spec(A)).$$
Since $\B$ is hereditary, $A(U) \in \B$ for all $U \in \CC_X$ and so $A$ is
locally $\B$-valued.
\qed \enddemo

%%%%%%%%%%%%%%%%%%%%%%%%%%%%%%%%%%%%%%%%%%%%%%%%%%%%%%%%%%%%%%%%%%%%%%%%%%%%%%%%
%%%%%%%%%%%%%%%%%%%%%%%%%%%%%%%%%%%%%%%%%%
\head 4. Asymptotic Riesz Representation Theorems \endhead
%%%%%%%%%%%%%%%%%%%%%%%%%%%%%%%%%%%%%%%%%%%%%%%%%%%%%%%%%%%%%%%%%%%%%%%%%%%%%%%%
%%%%%%%%%%%%%%%%%%%%%%%%%%%%%%%%%%%%%%%%%%

Throughout this section, we let $X$ denote a locally compact Hausdorff space
with Borel $\sigma$-algebra $\Sigma_X$. Let
$\CC_X \subset \Sigma_X$ denote the collection of all pre-compact open subsets
of $X$. And we let $\B \subset \B(\H)$
denote a hereditary $*$-subalgebra of the bounded operators on a fixed Hilbert
space $\H$.

\proclaim{Lemma 4.1} There is a bijective correspondence between
locally
$\B$-valued Borel POVMs $A : (\Sigma_X,\CC_X) \to  (\B(\H),\B )$
and positive linear maps $Q : C_0(X) \to \B$. This correspondence
is given
by $$Q(f) = \int_X f(x) \, dA(x) .\tag 4.1.1$$
\endproclaim

\demo{Proof} In view of Theorem A.3 we only need to check that the locally
$\B$-valued condition 
corresponds to
$Q(C_0(X)) \subset \B$. Suppose $A(\CC_X) \subset \B$. Let $f \in C_c(X)$ be
compactly supported.
Since $Q$ is positive linear, it suffices to assume $f \geq 0$. Let $K =
\supp(f)$ which is a compact subset
of $X$. By local compactness, there is an open subset $U \in \CC_X$ such that $K
\subset U$. By the Extreme
Value Theorem there is a $C> 0$ such that $0 \leq f \leq C \chi_U$. Since $Q$ is
positive,
$$0 \leq Q(f) = \int_X f \, dA \leq C \int_X \chi_U \, dA = C A(U) \in \B$$
by hypothesis. Since $\B$ is hereditary, $Q(f) \in \B$. 

Conversely, suppose $Q : C_0(X) \to \B$ is positive linear and given by formula
(4.1.1). Then $Q$ defines 
a positive map $Q : ( B_b(X),  C_0(X)) \to  (\B(\H), \B)$. Let $U \in
\CC_X$ be a pre-compact
open subset. Since $X$ is completely regular, we have by Urysohn's  Lemma a
continuous function $f \in C_c(X)$ with $0 \leq f \leq 1$
such that $\overline{U} = f^{-1}(1)$. Thus, $ 0 \leq \chi_U \leq f$
and so
$0 \leq A(U) = Q(\chi_U) \leq Q(f) \in \B$. Thus, $A(U) \in \B$ and so $A$ is
locally
$\B$-valued as desired.
 \qed \enddemo

Define $B_0(X)$ to be the $C^*$-subalgebra of $B_b(X)$ generated by $\{\chi_U :
U \in \CC_X\}$.
If $X$ is also $\sigma$-compact, a paracompactness argument then shows that
$C_0(X) \subset B_0(X)$ as a closed
(but not necessarily hereditary) $*$-subalgebra. (Recall that if $f \in C_c(X)$
is compactly supported, then
$\Interior(\supp(f)) \in \CC_X$.) The following is our main result.

\proclaim{Theorem 4.2} If $X$ is $\sigma$-compact, there is a bijective
correspondence between
positive asymptotic morphisms
$$\{Q_\hbar\} :( B_0(X),  C_0(X)) \to  (\B(\H), \B)$$
and locally $\B$-valued Borel asymptotic spectral measures
$$\{\Ah\} : ( \Sigma_X, \CC_X ) \to (\B(\H), \B).$$
This correspondence is given by
$$Q_\hbar(f) = \int_X f(x) \, d\Ah(x) \tag{4.1}$$
\endproclaim

\demo{Proof} Let $\{Q_\hbar\} : (B_0(X),  C_0(X)) \to  (\B(\H), \B)$ be
a positive asymptotic morphism. By the lemma, there is a locally $\B$-valued
family of 
POVMs $\{\Ah\} : ( \Sigma_X, \CC_X ) \to (\B(\H), \B)$ such that (4.1)
holds for all $\hbar > 0$ and $f \in B_0(X)$. For each $U \in \CC_X$ we
have that $\hbar \mapsto A_\hbar(U) = Q_\hbar(\chi_U)$ is continuous
by condition (2.1.b.). Also, we have that for any $U \neq \emptyset \in \CC_X$
$$\limsup_{\hbar \to 0} \|\Ah(U)\| = \limsup_{\hbar \to 0} \|Q_\hbar(\chi_U)\|
\leq
\|\chi_U\| =1.$$
Since $X$ is $\sigma$-compact, there is an increasing sequence $\{U_n\} \subset
\CC_X$ of pre-compact open 
subsets such that $X = \cup_1^\infty U_n$. By Theorem 18 \cite{Be}, for all
$\hbar > 0$,
$\Ah$ is a ``regular'' Borel POVM, so
$$\Ah(X) = \LUB \{\Ah(U_n) : n \in \N\}$$
(in the sense of positive operators). It follows that $\limsup_{\hbar \to
0}\|\Ah(X)\| \leq 1$ as desired.
Now let $U_1, U_2 \in \CC_X$. Since characteristic functions satisfy 
$\chi_{U_1 \cap U_2} = \chi_{U_1} \chi_{U_2}$
we then have by asymptotic multiplicativity (2.1.c) that
$$\lim_{\hbar \to 0} \|\Ah(U_1 \cap U_2) -
\Ah(U_1)\Ah(U_2)\| = 
\lim_{\hbar \to 0} \|Q_\hbar(\chi_{U_1}\chi_{U_2}) - 
Q_\hbar(\chi_{U_1})Q_\hbar(\chi_{U_2})|| = 0.$$
Thus, the family $\{\Ah\}$ is a locally $\B$-valued ASM on $X$.

Conversely, let $\{\Ah\} :( \Sigma_X, \CC_X ) \to (\B(\H), \B)$ be a
locally $\B$-valued ASM on $X$. 
Define the family of maps $\{Q_\hbar\} : B_0(X) \to \B(\H)$ by equation (4.1).
Hence, each $Q_\hbar$ is positive linear and 
$$Q_\hbar : ( B_0(X),  C_0(X)) \to  (\B(\H), \B)$$
by Lemma 4.1. Let $S_0(X)$ denote the dense
$*$-subalgebra of $B_0(X)$ 
consisting of simple functions $f= \sum_{i = 1}^n a_i \chi_{U_i}$, where $U_i
\in \CC_X$. 
Asymptotic projectivity (3.1.c) and the calculation above then show
that for any simple functions $f, g \in S_0(X)$ we have
$$\lim_{\hbar \to 0} \|Q_\hbar(fg) - Q_\hbar(f) Q_\hbar(g) \| = 0.$$
Also, for any such simple function $f \in S_0(X)$,
$$\hbar \to Q_\hbar(f) = \sum_1^n a_i \Ah(U_i)$$
is continuous from $(0,1] \to \B$ by (3.1.b).
To conclude that $\{Q_\hbar\}$ is asymptotically multiplicative on the closure
$B_0(X)$ we need to 
show that it is bounded. By (A.3.1) we have that for any $f \in
B_0(X)$,
$$\left\| Q_\hbar(f)\right\| =  \left\| \int_X f(x) \, d\Ah(x)  \right\| \leq 2 \left\|f\right\|
\left\|A_\hbar(X)\right\|.$$
By condition (3.1.a) we then have
$$\limsup_{\hbar \to 0} \|Q_\hbar(f) \| \leq 2 \|f\| \limsup_{\hbar \to 0}
\|\Ah(X)\| \leq 2 \|f\|.$$
The result now follows since every bounded asymptotic morphism on a dense
\newline
$*$-subalgebra extends to the 
closure. \qed \enddemo

\proclaim{Corollary 4.3} Under the above hypotheses, equivalent Borel asymptotic
spectral measures correspond to equivalent positive asymptotic morphisms. Thus,
there is
a well-defined map $\(X, \B\) \to \lcl C_0(X), \overline{\B} \rcl_{acp}$ which
maps $\(\Ah\) \to \lcl Q_\h \rcl_{acp}$.
\endproclaim

\demo{Proof} Follows from the fact that $\Ah(U) = Q_\hbar(\chi_U)$ and any
two asymptotic morphisms equivalent on a dense subalgebra, are equivalent.
Also, since $C_0(X)$ is nuclear, the second statement follows from Lemmas 2.3
and 2.5.
\qed \enddemo

Let $C_\delta(X)$ denote a unital $C^*$-subalgebra of $C_b(X)$ such that $C_0(X)
\ideal C_\delta(X)$.
By the Gelfand-Naimark Theorem \cite{GN}, $C_\delta(X) \cong C(\delta X)$ for
some `continuous' compactification $\delta X \supseteq X$.

\proclaim{Corollary 4.4} Let $\I \ideal \B(\H)$ be an ideal. Every locally
$\I$-valued
full Borel asymptotic spectral measure $\{\Ah\}$
on $X$ determines a canonical relative asymptotic morphism (in the sense of
Guentner \cite{G2})
$$\{Q_\hbar\} : (C_\delta(X), C_0(X)) \to (\B(\H), \I) $$
for any continuous compactification $\delta
X$ of $X$.
\endproclaim

\definition{Definition 4.5} A family $\{\Ah\}_{\hbar > 0} : \Sigma \to \B(\H)$
of Borel POV-measures on $X$ will be called a
{\it $C_\delta$-asymptotic spectral measure} if the family of maps $\{Q_\hbar\}$
defined by equation (4.1)
determines an asymptotic morphism $\{Q_\hbar\} : C_\delta(X) \to \B(\H)$.
\enddefinition

The following proposition is then easy to prove using Theorem A.3 and the
results above.

\proclaim{Proposition 4.6} There is a one-one correspondence between
locally $\B$-valued \newline $C_\delta$-asymptotic spectral measures
$$\{\Ah\} : ( \Sigma_X, \CC_X ) \to (\B(\H), \B)$$ and positive
asymptotic morphisms $$\{Q_\hbar\} :  ( C_\delta(X), C_0(X)) \to  
(\B(\H), \B).$$
\endproclaim

%%%%%%%%%%%%%%%%%%%%%%%%%%%%%%%%%%%%%%%%%%%%%%%%%%%%%%%%%%%%%%%%%%%%%%%%%%%%%%%%
%%%%%%%%%%%%%%%%%%%%%%%%%%%%%%%%%%%%%%%%%%
\head 5. Examples and Applications \endhead
%%%%%%%%%%%%%%%%%%%%%%%%%%%%%%%%%%%%%%%%%%%%%%%%%%%%%%%%%%%%%%%%%%%%%%%%%%%%%%%%
%%%%%%%%%%%%%%%%%%%%%%%%%%%%%%%%%%%%%%%%%%

%%%%%%%%%%%%%%%%%%%%%%%%%%%%%%%%%%%%%%%%%%%%%%%%%%%%%%%%%%%%%%%%%%%%%%%%%%%%%%%%
%%%%%%%%%%%%%%%%%%%%%%%%%%%%%%%%%%%%%%%%%%
\subhead 5.1 Constructing ASMs via Quantum Noise Models \endsubhead
%%%%%%%%%%%%%%%%%%%%%%%%%%%%%%%%%%%%%%%%%%%%%%%%%%%%%%%%%%%%%%%%%%%%%%%%%%%%%%%%
%%%%%%%%%%%%%%%%%%%%%%%%%%%%%%%%%%%%%%%%%%

We give a general method for constructing asymptotic spectral measures from
spectral measures (on a possibly different measure space)
by adapting a convolution technique used to model noise and uncertainty in
quantum measuring devices. 
See Section II.2.3 of Busch et al \cite{BGL} for the relevant background
material.

Let $(X_1, \Sigma_1)$ and $(X_2, \Sigma_2)$ be measure spaces. Let $\E_2 \subset
\Sigma_2$.
Consider a family of maps
$$\{p_\hbar\} : \Sigma_2 \times X_1 \to [0,1]$$
such that the following conditions hold:
\roster
\item"a.)" For every $\omega \in X_1$, $\Delta \mapsto p_\hbar(\Delta, \omega)$
is a probability measure on $X_2$;
\item"b.)" For each $\Delta \in \E_2$, the map $\hbar \to p_\hbar(\Delta,
\cdot)$ is continuous $[0,1) \to B_b(X_1)$; 
\item"c.)" For every $\Delta_1, \Delta_2 \in \E_2$, $$\lim_{\hbar \to 0} \|
p_\hbar(\Delta_1, \cdot) p_\hbar(\Delta_2, \cdot) - p_\hbar(\Delta_1 \cap
\Delta_2, \cdot)
\|_\infty = 0$$ where $\| \cdot \|_\infty$ denotes the sup-norm on $B_b(X_1)$.
\endroster

Let $E : \Sigma_1 \to \B(\H)$ be a spectral measure on $X_1$. Define a family of
maps $\{\Ah\} : \Sigma_2 \to \B(\H)$ by the formula
$$\Ah(\Delta) = \int_{X_1} p_\hbar(\Delta, \omega) \, dE(\omega)$$
for any $\Delta \in \Sigma_2$.

\proclaim{Theorem 5.1.2} The family $\{\Ah\} : \Sigma_2 \to \B(\H)$ defines an
ASM on $(X_2, \Sigma_2, \E_2)$. If $E$ is normalized then $\{\Ah\}$ is
also normalized.
\endproclaim

\demo{Proof} The fact that each $\Ah$ is a POVM on $X_2$ is easy. Continuity in
$\hbar$ follows from condition $(b.)$ and the following estimate for $\Delta \in
\E_2$,
$$ \|\Ah(\Delta) - A_{\hbar_0}(\Delta) \| = \left \| \int_{X_1} (p_\hbar(\Delta,
\omega) - p_{\hbar_0}(\Delta, \omega)) \, dE(\omega) \right \| 
\leq \| p_\hbar(\Delta, \cdot) - p_{\hbar_0}(\Delta, \cdot)\|_\infty. $$
Now we need to prove asymptotic projectivity. Let $\Delta_1, \Delta_2 \in
\E_2$. Consider the calculation
$$\align
\| & \Ah(\Delta_1) \Ah(\Delta_2)  - \Ah(\Delta_1 \cap \Delta_2)\| =  \\
& = \left \| \int_{X_1} p_\hbar(\Delta_1, \omega) \, dE(\omega) \int_{X_1}
p_\hbar(\Delta_2, \omega) \, dE(\omega)  -
\int_{X_1} p_\hbar(\Delta_1 \cap \Delta_2, \omega) \, dE(\omega) \right \| \\
& = \left \| \int_{X_1} p_\hbar(\Delta_1, \omega) p_\hbar(\Delta_2, \omega)\,
dE(\omega) - \int_{X_1} p_\hbar(\Delta_1 \cap \Delta_2, \omega) \, dE(\omega)
\right \|
\\
& \leq \|p_\hbar(\Delta_1, \cdot) p_\hbar(\Delta_2, \cdot) - p_\hbar(\Delta_1
\cap \Delta_2, \cdot)\|_\infty \to 0
\endalign $$ as $\hbar \to 0$ by (c.).
We finish by showing that the mild normalization condition holds:
$$\Ah(X_2) = \int_{X_1} p_\hbar(X_2, \omega) \, dE(\omega) = \int_{X_1} 1 \,
dE(\omega) = E(X_1) \leq I$$
by condition (a.) above and the fact that $E(X_1)$ is a projection.
\qed \enddemo

Note that the inequalities in the previous proof require that $E$ be a PVM. (See
Theorems 15 and 16 \cite{Be} and Theorem A.4.) See Example 5.3 below for a 
concrete example of this smearing technique.

The physical interpretation (for finite systems) is that $p_\hbar$ models the
noise or uncertainty in interpreting the readings of a measurement.
For example, if $E$ has an eigenstate $\phi = E(\{\omega\}) \phi$, then the
expectation value of $\Ah(\Delta)$ when the system
is in state $\phi$ is given by
$$\< \phi | \Ah(\Delta) | \phi \>  = p_\hbar(\Delta, \omega).$$
Thus, $p_\hbar$ determines a (conditional) confidence measure of the system.

%%%%%%%%%%%%%%%%%%%%%%%%%%%%%%%%%%%%%%%%%%%%%%%%%%%%%%%%%%%%%%%%%%%%%%%%%%%%%%%%
%%%%%%%%%%%%%%%%%%%%%%%%%%%%%%%%%%%%%%%%%%
\subhead 5.2 Quasiprojectors and Semiclassical Limits \endsubhead
%%%%%%%%%%%%%%%%%%%%%%%%%%%%%%%%%%%%%%%%%%%%%%%%%%%%%%%%%%%%%%%%%%%%%%%%%%%%%%%%
%%%%%%%%%%%%%%%%%%%%%%%%%%%%%%%%%%%%%%%%%%

In this example, we show that the theory of ASMs can be used to study
semiclassical limits.
The relevant background for the
material in this section can be found in Chapters 10 and 11 of Omnes book
\cite{O}.
We first need the following well-known result which is an easy consequence of
the functional calculus and spectral mapping theorem. 
(See also Lemma 5.1.6. \cite{WO}.) It gives a rigorous statement of the
procedure used to ``straighten out'' quasiprojectors
into projections.

\proclaim{Lemma 5.2.1} Let  $\{a_\hbar: \hbar >0\}$ be a continuous family of
elements in a $C^*$-algebra $\B$
such that $0 \leq a_\hbar \leq 1$ for each $\hbar > 0$ and $$\lim_{\hbar \to 0}
\|a_\hbar - a_\hbar^2\| = 0.$$
There is a continuous family of projections $\hbar \mapsto e_\hbar = e_\hbar^* =
e_\hbar^2$ such that
$$\lim_{\hbar \to 0}\|a_\hbar - e_\hbar\| = 0.$$
\endproclaim

Let $(X, \Sigma, \E)$ be an asymptotic measure space.

\proclaim{Proposition 5.2.2} Let $\{\Ah\}$ be a normalized ASM on
$(X,\Sigma,\E)$. For each
subset $\Delta \in \E$ there is a continuous family of projections
$\{E_\hbar(\Delta)\}$ such
that
$$\lim_{\hbar \to 0} \| \Ah(\Delta) - E_\hbar(\Delta) \| = 0.$$
Moreover if $\Delta_1$ and $\Delta_2$ are disjoint measurable sets in $\E$ then
$$\lim_{\hbar \to 0} \| E_\hbar(\Delta_1) E_\hbar(\Delta_2) \| = 0.$$
\endproclaim

\demo{Proof} For each $\Delta \in \E$ we have by monotonicity and
normalization that
$0 \leq \Ah(\Delta) \leq I$ for all $\hbar > 0$. Setting $\Delta = \Delta_1 =
\Delta_2$
in the asymptotic projectivity condition (3.1.c) we have that
$$\lim_{\hbar \to 0} \| A_\hbar(\Delta) - \Ah(\Delta)^2\| = 0.$$
Now invoke the previous lemma to get the continuous family $\{E_\hbar(\Delta)\}$
of projections.  
If $\Delta_1 \cap \Delta_2 = \emptyset$ then by condition (3.1.c) again, we have
that 
$$\lim_{\hbar \to 0}\|\Ah(\Delta_1) \Ah(\Delta_2)\| = 0.$$
A simple triangle inequality argument plus normalization then shows that 
$$\lim_{\hbar \to 0} \|E_\hbar(\Delta_1) E_\hbar(\Delta_2)\| = 0$$
as was desired. \qed \enddemo

The relation to semiclassical limits occurs when we take $X$ to be the locally
compact
phase space of a classical system and $\B = \B_1(\H)$ to be the algebra of
trace-class
operators.

\proclaim{Proposition 5.2.3} Let $\{\Ah\}$ be a Borel ASM on $X$ with locally
compact
trace.
Then for any subset $\Delta \in \CC_X$ we have
$$\lim_{\hbar \to 0} \tr(\Ah(\Delta) - \Ah(\Delta)^2) = 0$$
and there is a unique integer $N_\Delta \in \N$ such that
$$N_\Delta = \lim_{\hbar \to 0} \tr(\Ah(\Delta)).$$
Moreover, this integer is constant on the asymptotic equivalence class of
$\{\Ah\}$.
\endproclaim

\demo{Proof} The first limit follows from the continuity of
the trace.
Let $\{E_\hbar(\Delta)\}$ be the projections from the previous result. Since 
$$\Ah(\Delta) \in \B_1(\H) \subset \overline{\B_1(\H)} = \K(\H)$$ it follows
that $E_\hbar(\Delta) \in \K(\H)$,
i.e., $\hbar \mapsto E_\hbar(\Delta)$ is a continuous family of compact (hence,
finite rank)
projections. Therefore, since the rank of a projection is a continuous invariant
\cite{D}
$$\lim_{\hbar \to 0} \tr(\Ah(\Delta)) = \lim_{\hbar \to 0} \tr(E_\hbar(\Delta)) 
= \rank(E_{\hbar_0}(\Delta)) =_{\text{def}} N_\Delta$$
for any $\hbar_0 > 0$. The last statement follows again by continuity of the
trace.
 \qed \enddemo

Suppose $X$ denotes the position-momentum phase space $(x, p)$ of a
particle. Let $\{\Ah\}$ be a
locally compact trace Borel ASM on $X$. A bounded rectangle $R$ in phase space
with center
$(x_0, p_0)$ and
sides $2 \Delta x$ and $2 \Delta p$ can then be used to represent a classical
property asserting
the simultaneous existence of the position and momentum $(x_0, p_0)$ of the
particle with 
given error bounds $(\Delta x, \Delta p)$ on measurement. The nonnegative
integer $N_{R}$ which
satisfies $$N_{R} = \lim_{\hbar \to 0} \tr (\Ah(R))$$
can then be interpreted as the number of semiclassical states of the particle
bound in the rectangular box $R$,
which is familiar from  elementary statistical mechanics. We then have that
$$\aligned \tr(\Ah(R) - \Ah(R)^2 ) & = N_R O(\hbar) \\
          \tr(E_\hbar(R) - \Ah(R)) & = N_R O(\hbar).
\endaligned$$ Thus, $\hbar$ represents a {\it classicity parameter}. When $\hbar
\approx 0$ is small, the 
quantum representation of the classical property is essentially correct and when
$\hbar \approx 1$ the classical property
has essentially no meaning from the standpoint of quantum mechanics. Since these
relations are preserved 
on equivalence classes, ``a classical property corresponding to a sufficiently
large a priori bounds $\Delta x$
and $\Delta p$ is represented by a set of equivalent quantum projectors''
\cite{O},
i.e., equivalent locally compact trace ASMs. 
In addition, if $R_1$ and $R_2$ are disjoint rectangles, representing
distinct classical properties, then we have that
$$\| \Ah(R_1) \Ah(R_2) \| = O(\hbar)$$
and so ``two clearly distinct classical properties are (asymptotically) mutually
exclusive when considered as quantum properties'' \cite{O}.

%%%%%%%%%%%%%%%%%%%%%%%%%%%%%%%%%%%%%%%%%%%%%%%%%%%%%%%%%%%%%%%%%%%%%%%%%%%%%%%%
%%%%%%%%%%%%%%%%%%%%%%%%%%%%%%%%%%%%%%%%%%
\subhead 5.3 Unsharp Spin Measurements of Spin-$\frac12$ Systems \endsubhead
%%%%%%%%%%%%%%%%%%%%%%%%%%%%%%%%%%%%%%%%%%%%%%%%%%%%%%%%%%%%%%%%%%%%%%%%%%%%%%%%
%%%%%%%%%%%%%%%%%%%%%%%%%%%%%%%%%%%%%%%%%%

In this example, we give a geometric classification of certain asymptotic
spectral measures associated to pure spin-$\frac12$ particles. 

Recall that pure spin systems are represented by the Hilbert space $\H = \C^2$
\cite{BGL,S}. We then have $\B(\H) \cong M_2(\C)$. The Pauli spin operators
$\sigma_1, \sigma_2, \sigma_3$
are the $2 \times 2$ matrices 
$$\sigma_1 = \pmatrix 0 & 1 \\ 1 & 0 \endpmatrix, \sigma_2 = \pmatrix 0 & -i \\
i & 0 \endpmatrix, \sigma_3 = \pmatrix 1 & 0 \\ 0 & -1 \endpmatrix $$
which satisfy the relations:
\roster
\item"{$\bullet$}" $\sigma_i^* = \sigma_i, \sigma_i^2 = I$
\item"{$\bullet$}" $\sigma_i \sigma_j = - \sigma_i \sigma_j$ for $i \neq j$
\item"{$\bullet$}" $\sigma_i \sigma_j = i \epsilon_{ijk}\sigma_k$ for $i \neq j$
\endroster
where $I$ denotes the identity operator. A density  operator (or state) on $\H$
is a positive matrix $\rho \geq 0$ with trace one.
A fundamental result in the theory is the following.

\proclaim{Lemma 5.3.1} Any density operator $\rho$ on $\H$ can be written
uniquely in the
form
$$\rho = \rho(\vec x) = \frac12 (I + \vec x \cdot \vec \sigma), \quad \vec x =
(x_1, x_2, x_3)  \in
\R^3, \quad \|\vec x \| \leq 1,$$
where $\vec x \cdot \vec \sigma = x_1 \sigma_1 + x_2 \sigma_2 + x_3 \sigma_3$
and $\|\vec x\|^2 = x_1^2 + x_2^2 + x_3^2$.
Moreover, $\rho$ is a one-dimensional projection
iff $\vec x$ is a unit vector $\|\vec x \| =1$.
\endproclaim

\definition{Definition 5.3.3} A {\it spin} POVM on $X_2 = \{-\frac12,
+\frac12\}$ is a normalized POVM $A = \{A^+, A^-\}$ such that
$\tr(A^\pm) = 1$, where $A^\pm = A(\{\pm\frac12\})$. Thus $A^\pm \geq 0$ is a
density operator and $A^+ + A^- = I$. An asymptotic spectral measure $\{\Ah\}$
on $X_2$
will be called {\it spin} if each $A_\hbar$ is a spin POVM.
\enddefinition

Let $B^3 = \{\vec x \in \R^3 : \|\vec x\| \leq 1\}$ denote the closed unit
ball in $\R^3$. Let $S^2 = \partial B^3$ denote the unit sphere.
For each $\vec x \in B^3$ we obtain 
a spin POVM $A_x$ on $X_2$ by defining
$$A_x^\pm = \rho(\pm \vec x) = \frac12(I \pm \vec x \cdot \vec \sigma) \tag
5.3.4$$
which determines an ``unsharp'' spin observable.
Let $\lambda = \lambda(\vec x) = \|\vec x\|$ and define the quantities
$$r_x = \frac{1+ \lambda(\vec x)}{2} > \frac12,  \quad u_x =
\frac{1-\lambda(\vec x)}{2} < \frac12.$$
The quantity $r_x$ is called the degree of reality and $u_x$ is the degree of
unsharpness of the unsharp observable $A_x$ \cite{BGL,RK}.

\proclaim{Lemma 5.3.5} There is a bijective correspondence between spin POVMs $A
= \{A^+, A^-\}$ and points
$\vec x \in B^3$ in the closed unit ball of $\R^3$ given by (5.3.4).
Moreover, $A$ is a spectral measure if and only if  $\vec x \in S^2$ is a unit
vector.
\endproclaim

\demo{Proof} Follows from Lemma 5.3.1, Definition 5.3.3, and $A^- = I - A^+$.
\qed \enddemo

\proclaim{Theorem 5.3.6} There is a bijective correspondence between spin
asymptotic spectral measures $\{\Ah\} = \{\Ah^+, \Ah^-\}$ and 
continuous maps $\vec A : (0,1] \to B^3$ such that
$$\lim_{\hbar \to 0} \| \vec A(\hbar)\| = 1. \tag 5.3.6.1$$
This correspondence is given by the formula
$$\Ah^\pm = \frac12(I \pm \vec A(\hbar) \cdot \vec \sigma).$$ 
\endproclaim

\demo{Proof} By the Lemma we only need to prove continuity in $\hbar$ and that
asymptotic projectivity corresponds to condition (5.3.6.1)
above. By the properties of the Pauli spin operators above, we can show by
direct
computations (see also formulas (66a) and (66b) in \cite{S}) that
$$4[(A_\hbar^+)^2 - A_\hbar^+] + (1 - \|\vec A(\hbar)\|^2) I = 0$$
and
$$\|\vec A(\hbar) - \vec A(\hbar_0)\|^2 =  - \det((\vec A(\hbar)-\vec
A(\hbar_0)) \cdot \vec \sigma) = 4 \det(A_\hbar^+ -A_{\hbar_0}^+).$$
The result now easily follows.\qed \enddemo

Thus, we can geometrically realize the space of spin asymptotic spectral
measures as the space of continuous
paths in the closed unit ball of $\R^3$ which asymptotically approach the  unit
sphere, i.e. they are ``asymptotically sharp.'' Note that this provides
nontrivial examples of asymptotic spectral measures which do not converge to 
a fixed spectral measure.

Let $\vec n$ be a unit vector and define $\vec A(\hbar) = (1 -
\hbar)\vec n$. The associated spin asymptotic
spectral measure given by  
$$A_h^\pm = \frac12(I \pm (1-\hbar)\vec n \cdot \vec \sigma)$$
is used by Roy and Kar \cite{RK} to analyze eavesdropping strategies in
quantum cryptography using EPR pairs of correlated
spin-$\frac12$ particles. Violations of Bell's inequality occur when the
parameter $\hbar > 1 - \sqrt{2}(\sqrt{2}-1)^\frac12$.

This spin ASM is also obtained by the asymptotic smearing construction in 5.1.
Let
$E^\pm = A_0^\pm$ be the spectral measure
associated to the unit vector $\vec n$. Define the family $\{p_\hbar\} : \Cal
P(X_2) \times X_2
\to [0,1]$ by the formula $p_\hbar(\Delta, j) = \sum_{i \in \Delta}
\lambda^\hbar_{ij}$ where $(\lambda^\hbar_{ij})$ is the stochastic matrix
$$(\lambda^\hbar_{ij}) = \pmatrix  1 - \frac{\hbar}{2} & \frac{\hbar}{2} \\  
\frac{\hbar}{2} & 1 - \frac{\hbar}{2} \endpmatrix.$$
One can then verify that
$$\Ah^\pm = \sum_{l = \mp \frac12} p_\hbar(\{\pm \frac12\}, l) E^\mp.$$

\proclaim{Corollary 5.3.7} Two spin asymptotic spectral measures $\{\Ah\}$ and
$\{B_\hbar\}$ are equivalent if and only if their 
associated maps $\vec A, \vec B : (0, 1] \to B^3$ are asymptotic, i.e., 
$$\lim_{\hbar \to 0} \|\vec A(\hbar) - \vec B(\hbar)\| = 0.$$
\endproclaim

\demo{Proof} $\| \vec A(\hbar) - \vec B(\hbar) \|^2 = 4 \det(A_\hbar^+ -
B_\hbar^+)$. \qed \enddemo

%%%%%%%%%%%%%%%%%%%%%%%%%%%%%%%%%%%%%%%%%%%%%%%%%%%%%%%%%%%%%%%%%%%%%%%%%%%%%%%%
%%%%%%%%%%%%%%%%%%%%%%%%%%%%%%%%%%%%%%%%%%
\subhead 5.4 Strong Deformation Quantization \endsubhead
%%%%%%%%%%%%%%%%%%%%%%%%%%%%%%%%%%%%%%%%%%%%%%%%%%%%%%%%%%%%%%%%%%%%%%%%%%%%%%%%
%%%%%%%%%%%%%%%%%%%%%%%%%%%%%%%%%%%%%%%%%%

Let $X$ be a locally compact space. Let $\B$ be a $C^*$-algebra. A {\it strong
deformation} from $X$ to $\B$ is a continuous field \cite{D} of
$C^*$-algebras $\{\B_\hbar : \hbar \in [0,1]\}$ such that $\B_0 = C_0(X)$ and
$$\{\B_\hbar \, | \, \hbar > 0\} \cong \B \times (0,1].$$
Here we give a measure-theoretic criterion, based on $E$-theory arguments, for
when a locally $\B$-valued Borel ASM $\{\Ah\}$ on $X$ determines
a strong deformation from $X$ to $\B$, where $\B$ is a hereditary
$C^*$-subalgebra of $\B(\H)$. First, we make the following general definition.

\definition{Definition 5.4.1} Let $\{\Ah\}$ be an ASM on $(X, \Sigma, \E)$. We
will call
$\{\Ah\}$ {\it injective} if
$$\liminf_{\hbar \to 0} \| \Ah(\Delta) \| > 0 $$
for all nonempty subsets $\Delta \neq \emptyset$ in $\E$.
\enddefinition

Thus, if $\{\Ah\}$ is a locally $\B$-valued Borel ASM on $X$, then by local
compactness and monotonicity,
injectivity is equivalent to $$\liminf_{\hbar \to 0} \|\Ah(U)\| > 0$$
for all nonempty open subsets $U \neq \emptyset$ of $X$.
Let $\{Q_\hbar\} : C_0(X) \to \B$ be the
associated asymptotic morphism given by Theorem 4.2. 
Recall that $\{Q_\hbar\}$ is called injective \cite{L1} if 
$$\liminf_{\hbar \to 0} \| Q_\hbar(f) \| > 0$$
for all $f \neq 0$ in $C_0(X)$. By the results in \cite{CH,L1,DL}, (weakly)
injective
asymptotic morphisms determine strong deformations from $X$
to $\B$.

\proclaim{Theorem 5.4.2} Let $\{\Ah\}$ be an injective locally $\B$-valued Borel
ASM on
$X$. Then the associated asymptotic morphism
$\{Q_\hbar\} : C_0(X) \to \B$ is injective and so satisfies the continuity
condition
$$\|f\| = \lim_{\hbar \to 0} \|Q_\hbar(f)\|$$
for all $f \in C_0(X)$.  Hence, there is an associated strong deformation from
$X$ to $\B$.
\endproclaim

\demo{Proof} Let $f \neq 0$ be in $C_0(X)$. Thus, there is an $x_0 \in X$ such
that $|f(x_0)| > C > 0$. Since $\{Q_\hbar\}$ is positive linear, without
loss of generality, we may assume $f \geq 0$ and so $f(x_0) > C > 0$. Let $U \in
\CC_X$ be the pre-compact open subset
of $X$ defined by $U = \{x \in X : f(x) > C\}.$
Then $C \chi_U \leq f$ and so for all $\hbar > 0$ we have that
$$ C \Ah(U) = \int_X C \chi_U\, d\Ah \leq \int_X f \, d\Ah =  Q_\hbar(f)$$ which
implies that
$$0 < |C| \liminf_{\hbar \to 0} \| A_\hbar(U) \| \leq \liminf_{\hbar \to 0}
\|Q_\hbar(f) \|.$$
It follows that $\{Q_\hbar\} : C_0(X) \to \B$ is
injective and so by Lemma 3 \cite{L1} 
$$\|f\| = \lim_{\hbar \to 0} \|Q_\hbar(f)\|$$
for all $f \in C_0(X)$. Thus, $\{Q_\hbar\}$ is the asymptotic morphism
associated to a strong deformation from $X$ to $\B$. \qed \enddemo

The continuous sections of the field $\{\Ah\}$ are then determined by the
equivalence
class $\lcl Q_\hbar \rcl_{a(cp)}$ of the associated asymptotic morphism
$\{Q_\hbar\} : C_0(X) \to \B$.

%%%%%%%%%%%%%%%%%%%%%%%%%%%%%%%%%%%%%%%%%%%%%%%%%%%%%%%%%%%%%%%%%%%%%%%%%%%%%%%%
%%%%%%%%%%%%%%%%%%%%%%%%%%%%%%%%%%%%%%%%%%
\subhead 5.5 Wick Quantization on Bosonic Fock Space \endsubhead
%%%%%%%%%%%%%%%%%%%%%%%%%%%%%%%%%%%%%%%%%%%%%%%%%%%%%%%%%%%%%%%%%%%%%%%%%%%%%%%%
%%%%%%%%%%%%%%%%%%%%%%%%%%%%%%%%%%%%%%%%%%

The background material for this section can be found in Guentner\cite{G1}.
Let $\H = L^2(\C, d\mu(z))$ denote the Hilbert space of measurable
complex-valued functions on the complex plane $X
= \C$ which are square-integrable with respect to the normalized Gaussian
measure $d\mu(z) = \pi^{-1}e^{-|z|^2} d\lambda(z) = \pi^{-1}k(z, z)d\lambda(z),$
where $k(z,w) = e^{z\bar{w}}$ denotes the {\it Bergman kernel} and $d\lambda(z)$ denotes Lebesgue measure. The
(bosonic) {\it Fock space}  is the closed subspace $\F \subset \H$ consisting of analytic functions.
For any bounded Borel function $f \in B_b(\C)$, the {\it Wick operator}
$T_f : \F \to \F$ of $f$ is the integral operator
defined by
$$T_f(\phi) = \int_\C k(z,w) f(w) \phi(w) d\mu(z),$$
for all $\phi \in \F$.

\proclaim{Lemma 5.5.1} For each $f \in L^2(\C, d\lambda) \cap B_b(\C)$ the
operator $T_f \in \K(\F)$.
\endproclaim

\demo{Proof} Follows from  the calculations in the proof of Proposition 3.2
\cite{G1} 
\qed \enddemo

We define the Wick quantization map 
$$Q^W : C_b(\C) \to \B(\F) : f \mapsto Q^W(f) = T_f.$$ 
Let $P : \H \to \F$
denote the orthogonal projection. We can then define a POV-measure $A^W :
\Sigma_\C \to \B(\F)$ by
$$A^W(\Delta) = P \circ \chi_\Delta$$
where $\chi_\Delta$ denotes (the operator on $\H$ of multiplication by) the
characteristic function $\chi_\Delta$.
Note that it is the compression of the PVM $\Delta \mapsto E(\Delta) =
\chi_\Delta$.

\proclaim{Lemma 5.5.2} The POVM $A^W :  ( \Sigma_\C, \CC_\C) \to 
 (\B(\F), \K(\F))$ is normalized and locally compact-valued. The associated
positive
linear map is the Wick-Toeplitz quantization $$Q^W : (C_b(\C), C_0(\C)) \to
(\B(\F), \K(\F)).$$
\endproclaim

\demo{Proof} Follows from the fact that $A^W(\Delta) = P \circ \chi_\Delta =
T_{\chi_\Delta} = Q^W(\chi_\Delta)$. 
When $U \in \CC_\C$ is pre-compact then $\chi_U \in  L^2(\C, d\lambda)
\cap B_b(\C)$ and so $A^W(U) = T_{\chi_U} \in \K(\F)$. Normalization
follows from $A^W(\C) = P = I_\F$. \qed 
\enddemo

For each $\hbar >0$ and $f \in B_b(\C)$ define
$$\alpha_\hbar(f)(z) = f(\hbar z)$$
for all $z \in \C$. We can then define a family of positive linear maps 
$$\{Q^W_\hbar\} : (C_b(\C), C_0(\C)) \to (\B(\F), \K(\F)) : f \mapsto
Q^W_\hbar(f) =
Q(\alpha_\hbar(f)).$$
Guentner \cite{G1} realized that to obtain an asymptotic morphism from the Wick
quantization we need to pass to a unital subalgebra of $C_b(\C)$
that still contains $C_0(\C)$ as an ideal. 

Let $\delta\C$ denote the compactification of the complex plane $\C$ by the
circle at infinity. The continuous functions
on $\delta\C$ are ``flat at infinity'' when restricted to $\C$. Let
$C_\delta(\C) = C(\delta\C)$.
We then have that $C_0(\C) \ideal C_\delta(\C) \subset C_b(\C)$. The following
result is a consequence of Proposition 4.6 above and Propositions 3.2 and 3.3
\cite{G1}.

\proclaim{Proposition 5.5.4} The family $\{Q^W_\hbar\}$ defines a
relative positive asymptotic morphism
$$\{Q^W_\hbar\} : (C_\delta(\C), C_0(\C)) \to (\B(\F), \K(\F))$$
whose associated $C_\delta$-asymptotic spectral measure $\{\Ah^W\}$ is
given by
$$\Ah^W(\Delta) = A^W(\hbar^{-1} \Delta) = P \circ
\alpha_\hbar({\chi_\Delta}).$$
\endproclaim

The restricted asymptotic morphism $\{Q^W_\hbar\} : C_0(\C) \to  \K(\F)$
determines an $E$-theory class, 
$$\lcl A^W_\hbar \rcl =_{\text{def}} \lcl Q^W_\hbar\rcl \in E(C_0(C), \K) \cong
E(C_0(\C), \C),$$
where we have used the matrix-stability of $E$-theory.

Let $\dbar = \frac12(\frac{\partial}{\partial x} +
\sqrt{-1}\frac{\partial}{\partial y})$ be the $\dbar$-operator
on $\C \cong \R^2$, considered as an unbounded elliptic differential operator
on the Hilbert space $H = L^2(\C, d\lambda(z))$. The formal adjoint of $\dbar$
is the operator
$-\partial$ where $\partial = \frac12(\frac{\partial}{\partial x} -
\sqrt{-1}\frac{\partial}{\partial y})$. It follows
that the $2 \times 2$ matrix operator
$$D = \pmatrix 0 & -\partial \\ \dbar & 0 \endpmatrix$$
determines a symmetric unbounded operator on $H \oplus H$ with bounded
propagation, and so is (essentially) self-adjoint.

By the results of Guentner \cite{G1,G2} the operator $D$ determines an
$E$-theory class denoted $\lcl \dbar \rcl \in
E(C_0(\C),
\C)$, which is the homotopy class of the asymptotic morphism determined by the
formula
$$C_0(\R) \otimes C_0(\C) \to C_0(\R) \otimes \K(H \oplus H) : f \otimes \phi
\mapsto M_\phi \circ f(\hbar D +
x\epsilon)$$ where $x \in \R$ and $\epsilon$ is the grading operator of the
$\Z_2$-graded Hilbert space $H \oplus H$.
A direct consequence of Proposition 5.5.4 above, Theorem 4.5 \cite{G1}, and the
excision property of relative $E$-theory
\cite{G2} is that the
$E$-theory classes of the Wick ASM above and the $\dbar$-operator are in fact
equal.

\proclaim{Theorem 5.5.5} $\lcl A^W_\hbar \rcl  =  \lcl \dbar \rcl \in E(C_0(\C),
\C) \cong \Z$.
\endproclaim

%%%%%%%%%%%%%%%%%%%%%%%%%%%%%%%%%%%%%%%%%%%%%%%%%%%%%%%%%%%%%%%%%%%%%%%%%%%%%%%%
%%%%%%%%%%%%%%%%%%%%%%%%%%%%%%%%%%%%%%%%%%
\head Appendix: POV-Measures and Quantum Mechanics \endhead
%%%%%%%%%%%%%%%%%%%%%%%%%%%%%%%%%%%%%%%%%%%%%%%%%%%%%%%%%%%%%%%%%%%%%%%%%%%%%%%%
%%%%%%%%%%%%%%%%%%%%%%%%%%%%%%%%%%%%%%%%%%

Let $X$ be a set equipped with a $\sigma$-algebra $\Sigma$ of subsets of $X$.
Let $\H$ be a separable Hilbert space with inner product $\< \cdot, \cdot \>$. A
{\it positive operator-valued measure} 
(POVM) on the measurable space $(X, \Sigma)$ is a mapping $A : \Sigma \to
\B(\H)$
which satisfies the following properties:

\roster
\item"$\bullet$" $A(\emptyset) = 0$
\item"$\bullet$" $A(\Delta) \geq 0$ for all $\Delta \in \Sigma$
\item"$\bullet$" $A(\sqcup_1^\infty \Delta _n) = \sum_1^\infty A(\Delta _n)$ for
disjoint measurable sets $\{\Delta _n\}_1^\infty \subset \Sigma$,
\endroster
where the sum converges in the weak operator topology \cite{Be,BGL,S}. Note that
$$0 \leq A(\Delta) \leq A(X) \leq \|A(X)\| < \infty$$ for all $\Delta \in
\Sigma.$
We will say that $A$ is {\it normalized} if  $A(X) = I_\H$. If each $A(\Delta )$
is
a projection in $\B(\H)$, i.e.,  $A(\Delta )^2 = A(\Delta)^* = A(\Delta ),$
then we call $A$ a {\it projection-valued measure} (PVM or {\it spectral
measure}) on $X$. This is equivalent to the condition that:
For all $ \Delta_1, \Delta_2 \in \Sigma$,
$$A(\Delta_1 \cap \Delta_2) = A(\Delta_1)A(\Delta_2) \tag A.1$$
See Berberian \cite{Be} for the basic integration theory of POVMs and Brandt
\cite{Br} for a short history of POVMs in quantum theory and an application
to photonic qubits in quantum information processing. The monographs
\cite{S,BGL}
give a thorough exposition of POVMs in foundational and operational
aspects of quantum physics.

Let $X$ be a locally compact Hausdorff topological space. Let $C_0(X)$ denote
the $C^*$-algebra of all continuous complex-valued
functions on $X$ which vanish at infinity. 

\definition{Definition A.2}
A {\bf general quantization} of $X$ on a Hilbert space $\H$ is a positive linear
map $Q : C_0(X) \to
\B(\H)$. If $X$ is compact, we require that $Q(1_X) = I_{\H}$. If
$X$ is a non-compact space, we require that $Q$ should have a unital extension
$Q^+ : C_0(X)^+
= C(X^+) \to \B(\H)$ which is a positive
linear map.
\enddefinition

Another reason for the importance of these operator-valued measures in
quantization is the following
generalized Riesz representation theorem for the dual of $C_0(X)$ (Compare
Proposition 1.4.8 \cite{L} and Theorem 19 \cite{Be}):

\proclaim{Theorem A.3} Let $\Sigma_X$ be the Borel $\sigma$-algebra on the space
$X$. There is a one-one
correspondence between positive linear maps  $Q : C_0(X) \to \B(\H)$ and
POV-measures
$A : \Sigma_X \to \B(\H)$, given by
$$Q(f) = \int_X f(x) \ dA(x).$$ 
The map $Q$ is a general quantization if and only if $A$ is a normalized POVM.
Moreover, $Q$ is a $*$-homomorphism if and only if $A$ is a spectral measure
(PVM).
\endproclaim 
\noindent The above integral is to be interpreted in the weak sense: For all $v,
w \in \H$,
$$\< Q(f)v , w\rangle = \int_X f(x) \< dA(x) v, w \>.$$
The map $Q$ then extends to $Q: B_b(X) \to \B(\H)$ and satisfies (Theorem 10
\cite{Be}): For all $f \in
C_0(X) \subset B_b(X)$,
$$\|Q(f)\| = \left \| \int_X f(x) \, dA(x) \right \| \leq 2 \|f\| \|A(X)\|. \tag
A.3.1$$
Thus, spectral measures (PVM's) correspond to representations of abelian
$C^*$-algebras
on Hilbert space. The fundamental result in the von Neumann formulation 
of quantum theory is the following Spectral Theorem of Hilbert.

\proclaim{Spectral Theorem A.4} 
Let $X = \R$. There is a one-one correspondence between Borel spectral measures
$A$ on $\R$ and
self-adjoint operators $T$ on the associated Hilbert space. This correspondence
is given by the formulas:
$$T = \int_{-\infty}^\infty \lambda \, dA(\lambda), \quad A(\Delta) =
\chi_\Delta(T),$$
where $\chi_\Delta$ denotes the characteristic function of the Borel set $\Delta
\subset \R$. 
\endproclaim

Let $A$ be a normalized POV-measure on the phase space $X$ of a quantum system.
The physical interpretation of the map $\Delta \mapsto
A(\Delta)$ is the probability that the physical system,  in a state represented
by a density operator $\rho$, is
localized in the subset $\Delta$ of the phase space $X$ is given by the number 
$$ P_\rho(\Delta) = \tr(\rho \circ A(\Delta)) = \tr \Big( \int_\Delta \rho \ dA
\Big).$$
The mean or vacuum expectation value of a quantum observable $T$ is then
computed by the formula
$$\langle T  \rangle = \tr(\rho T) = \tr \Big( \int_{-\infty}^\infty \lambda
\rho(\lambda) \ dA(\lambda) \Big),$$
where $\rho$ is the (normalized) probability density operator of the physical
system.

Note that according to the Naimark Extension Theorem \cite{RS}, every
POVM $A$ is the compression of a PVM $E$ defined on a
minimal extension $\H' \supset \H$. That is, $A(\Delta) = P E(\Delta) P$, where
$P : \H' \to \H$ is the orthogonal
projection. One could then try to compute the integrals
$\int_X f(x) \ dA(x)$ by computing $\int_X f(x) \ dE(x)$ on $\H$ and then
projecting back down to $\H$. There are two
problems with this \cite{S}. The first is that $\H'$ could have no physical
meaning,
thus making the analysis unsatisfying  to
the physicist. Also,  the integration process may not commute with the
projection process (e.g., when the associated operator is unbounded).

\vfill

%\eject

%%%%%%%%%%%%%%%%%%%%%%%%%%%%%%%%%%%%%%%%%%%%%%%%%%%%%%%%%%%%%%%%%%%%%%%%%%%%%%%%

\enddocument
\end